\documentclass[12pt,english]{article}

\usepackage[latin1]{inputenc}
\usepackage[english]{babel}
\usepackage[T1]{fontenc}
\usepackage{amsfonts}
\usepackage{latexsym}
\usepackage{amsthm}
\usepackage{amsmath}
\usepackage{mathrsfs}
\usepackage{amssymb}
\usepackage{graphpap}
\usepackage{graphicx}
\usepackage{epsfig}

\let\phi=\varphi
\def\Z{{\mathbb Z}}
\def\N{{\mathbb N}}
\def\eps{\varepsilon}
\def\IP{{\mathbb P}}
\def\IE{{\mathbb E}}

\def\G{{\mathcal G}}
\def\V{{\mathcal V}}
\def\Ed{{\mathcal E}}

\def\L{{\mathcal L}}
\def\T{\mathcal{T}}
\def\Po{{\mathtt P}_{\omega}}
\def\Eo{{\mathtt E}_{\omega}}

\def\P{\textbf{P}}
\def\E{\textbf{E}}

\def\v{{\bf v}}
\def\u{{\bf u}}
\def\w{{\bf w}}
\def\x{{\bf x}}
\def\y{{\bf y}}
\def\z{{\bf z}}
\def\r{{\bf r}}
\def\w{{\bf w}}

\newcommand{\argmin}{\mathop{\mathrm{arg\,min}}}

\newtheorem{theo}{Theorem}[section]

\newtheorem{rem}{Remark}[section]
\newtheorem{df}{Definition}[section]
\newtheorem{prop}{Proposition}[section]

\allowdisplaybreaks

\title{Spiders in random environment}
\author{Christophe Gallesco$^{~1}$, Sebastian M\"uller$^{~2}$, Serguei Popov$^{~3}$,\\ Marina Vachkovskaia$^{~4}$}

\begin{document}

\bibliographystyle{plain}

\maketitle
{\footnotesize 

\noindent $^{~1}$ Instituto de Matem\'atica e Estat\'istica, Universidade de
S\~ao Paulo, rua do Mat\~ao 1010, CEP 05508--090, S\~ao Paulo, SP, Brazil

\noindent e-mail: \texttt{gallesco@ime.usp.br}

\noindent $^{~2}$  L.A.T.P. / C.M.I., Universit\'e de Provence
39 rue Joliot Curie, 13453 Marseille cedex 13, France 

\noindent e-mail: \texttt{mueller@cmi.univ-mrs.fr}\\
url: \texttt{http://www.latp.univ-mrs.fr/~muller}

\noindent $^{~3,4}$  Department of Statistics, Institute of Mathematics, Statistics and Scientific Computation, University of Campinas--UNICAMP, rua S\'ergio Buarque de Holanda 651, 13083--859, Campinas SP, Brazil

\noindent e-mail: \texttt{popov@ime.unicamp.br}; \texttt{marinav@ime.unicamp.br}\\
url: \texttt{http://www.ime.unicamp.br/$\sim$popov}; \texttt{http://www.ime.unicamp.br/$\sim$marinav}

}

\begin{abstract}
A spider consists of several, say $N$, particles. Particles can jump
independently according to a random walk if the movement does not
violate some given restriction rules. If the movement violates a rule
it is not carried out. We consider random walk in random environment
(RWRE) on $\Z$ as underlying random walk. We suppose the environment
$\omega=(\omega_x)_{x \in \Z}$ to be elliptic, with positive drift and nestling, so
that there exists a unique positive constant $\kappa$ such that
$\E[((1-\omega_0)/\omega_0)^{\kappa}]=1$. The restriction rules are kept very general; we only
assume transitivity and irreducibility of the spider. The main
result is that the speed of a spider is positive if $\kappa/N>1$ and null if $\kappa/N<1$. In particular, if $\kappa/N <1$ a spider has null speed but the
speed of a (single) RWRE is positive.
 \end{abstract}

\section{Introduction and results}
To begin with, let us give a simple example of  a spider. Imagine there are two particles performing nearest neighbor
random walks on $\Z$ in continuous time. These particles are tied together with a rope of a certain length $s\in \N$. As long as the rope is not tight their movements are independent. If the rope is tight (the two particles
are at  a distance $s$ from each other) the rope prevents the particles to jump {\it away} from each other.

In these notes, we consider a spider on $\Z$ in a random environment. First, suppose that
$\omega:=(\omega_x)_{x\in \Z}$  is a sequence of positive i.i.d.\ random variables taking values in $(0,1)$. We
denote by $\P$ the distribution of $\omega$ and by $\E$ the corresponding expectation. In the
example above, we first choose an environment $\omega$ at random according to the law $\P$ and we describe the
position of our two particles by the vector $S(t)=(S_1(t), S_2(t))$ where $S_i(t)$, $i=1,2$, is the position of
particle $i$ at time $t$. As long as $|S_{1}(t)-S_{2}(t)|<s$, the two particles behave like two independent random walks in random environment. If $|S_{1}(t)-S_{2}(t)|=s$, their movements are dependent in order to prevent that $|S_{1}(t)-S_{2}(t)|>s$. For instance, let the first particle be in $x_1$ and the second in $x_2$.
Then, if $|x_1-x_2|<s$ the first particle jumps  to $x_1+1$ with rate $\omega_{x_1}^+:=\omega_{x_1+1}$ or to
site $x_1-1$ with rate $\omega_{x_1}^-:=1-\omega_{x_1}$. Analogously the second one moves to $x_2+1$ with rate
$\omega_{x_2}^+$ or to site $x_2-1$  with $\omega_{x_2}^-$. If $|x_1-x_2|=s$ and $x_1<x_2$ the first leg may
only jump to the right with rate $\omega_{x_1}^+$ and the second to the left with rate $\omega_{x_2}^-$. In the
case $x_1<x_2$ the roles of the two legs are interchanged.

More generally we can consider a spider with $N$ legs, that is to say $N$ interacting particles. The particles move independently as long as their movement does not violate some restriction rules concerning their positions. In
this case we denote by $S(t)=(S_1(t),S_2(t),\dots,S_N(t))$ the positions of the $N$ particles at time $t$ where $S_i(t)$ represents the position of particle $i$ at time $t$.

This model gained recently an interest in evolutionary dynamics and molecular cybernetics. At the moment, to our
knowledge, there are just a few theoretical papers on this model. In~\cite{AKM}, Antal, Krapivsky and Mallick
obtained the speed and diffusion constants for 1-dimensional spiders and in~\cite{AK} Antal and Krapivsky made the first study
for non-Markovian spiders. In~\cite{GMP}, Gallesco, M\"uller and Popov study qualitative properties, as recurrence, transience, ergodicity
and positive rate of escape of spiders in a quite general setting. We refer to the lecture notes of Zeitouni~\cite{Z} for a general overview
on random walks in random environments (RWRE).  The main result of this paper, Theorem~\ref{theo1}, is that in
random environment on $\Z$ the speed of a spider may be zero even if the speed of a (single) RWRE is positive.
This is in contrast with the results in~\cite{AKM} and in Section 4.1 of~\cite{GMP} that the positive speed of a
homogeneous random walk implies positive speed of the spider.

It is convenient to adopt the following notations of~\cite{GMP}. Recall that the spider is described through
$S(t)=(S_1(t),\ldots,S_N(t))$. The first leg defines the position of the spider: the position of the spider at
time $t$ is $S_1(t)$. The spider is defined through a set $L$ of local configurations at $0$, that is a finite
subset of $\{(x_1,x_2,\dots,x_N): x_1=0,~x_2,\dots,x_N \in \Z\}$. Actually, the set $L$ corresponds  to all
possible configurations for the spider  at position $0$. Since in this note we only consider transitive spiders
the set of local configurations at position $x$ (that is when $x_1=x$) can be obtained by translating the set
$L$ by $x$. Denoting by $\Theta_x$ the shift by $x$, we have
\[L_x=\Theta_xL=\{(x,\ x_2, \dots,x_N)\in \Z^N: (0, x_2-x,\dots,x_N-x)\in L\}.
\]
 Let
\[\V=\bigcup_{x\in \Z}L_x.
\]
For elements in $\V$ we write $\x=(x_1,\dots,x_N)$ and $\y=(y_1,\dots,y_N)$. The transition rates
$Q^S=Q^S(\omega)=(q^S(\x,\y))_{(\x,\y) \in \V^2}$ of the spider are defined as follows: let $\x,\y\in\V$ then
\begin{itemize}
\item if $\|\x-\y\|=1$ (where $\|\cdot\|$ is the usual $\ell_1$-norm) and $i$ is the coordinate such that
$x_i\neq y_i$,
\[
q^S(\x,\y)= \left\{
    \begin{array}{ll}
        \omega^+_{x_i} & \mbox{if}~y_i=x_i+1\\
        \omega^-_{x_i} & \mbox{if}~y_i=x_i-1,\\
    \end{array}
\right.
\]
\item otherwise,
\[
q^S(\x,\y)=0.
\]
\end{itemize}

Now, following \cite{GMP}, we define the spider graph. For a given realization
$\omega$ of our environment, define the graph $\G=\G(\omega)=(\V,\Ed(\omega))$ such that an edge $e=(\x,\y)\in
\V \times \V$ belongs to $\Ed(\omega)$ if and only if $q^S(\x,\y)>0$. As the sequence $\omega$ takes values in $(0,1)^{\Z}$, 
the spider graph is deterministic. In the rest of these notes we will assume the irreducibility  of the spider
walk which is implied by the two following conditions on $\G$ for almost all realizations of $\omega$:
\begin{itemize}
\item [(i)] $L$ is a connected subgraph of the spider graph $\G$, \item [(ii)] there exists at least one edge between $L$ and $L_1$.
\end{itemize}
Condition (i) is not necessary for the irreducibility of the spider, nevertheless it is assumed in this
stronger form to reduce the technical part of the proofs. We assume the following conditions on our random
environment:
\begin{itemize}
\item [(iii)] $\E[\ln \rho_0]<0$, with $\rho_0:=\frac{\omega_0^-}{\omega_0^+}$, \item [(iv)] there
exists~$0<\delta<1/2$ such that $\P[\delta \leq \omega^+_0 \leq 1-\delta]=1$, \item [(v)] $\P[\omega_0^+>1/2]>0$
and $\P[\omega_0^+\leq 1/2]>0$.
\end{itemize}
Condition (iii) implies that the RWRE is transient to the right, see Solomon \cite{solom}. Condition (iv) is the usual uniform ellipticity condition. Condition (v)
corresponds to the fact that our environment is nestling. Observe that for non-nestling random environments it
is possible to show that every spider, satisfying (i)+(ii), has positive speed.
Furthermore, conditions (iii)-(v) imply that there exists a unique $\kappa>0$, such that
\[
\E[\rho_0^{\kappa}] = 1.
\]

We denote by $\Po^{\x}$ the quenched law of the spider starting at $\x$ in the environment $\omega$ and by $\Eo^{\x}$ the
corresponding expectation. Finally, we denote by $\IP^{\x}:=\P\times \Po^{\x}$ and $\IE^{\x}$ the annealed
probability and expectation for the spider starting at $\x$.
\medskip

We define the speed of a spider as
\begin{equation*}
v=\lim_{t \to \infty} \frac{S_1(t)}{t}
\end{equation*}
if the limit exists.
Let us consider a spider starting at some initial position $\x_0\in L$ and define the stopping time
\[
\T:= \inf \{s>0: S_1(s)>0 \phantom{*}\mbox{and}\phantom{*}S(s)=\Theta_{S_1(s)} \x_0\}.
\]
The main result of these notes is the following theorem.

\begin{theo}
\label{theo1}  Consider a spider with $N$ legs. Under conditions (i)-(v), the speed $v$ of the spider is well-defined and we have $\IP$-a.s.
\[
 v=\frac{\IE[S_1(\T)]}{\IE[\T]}>0 \quad \mbox{if\phantom{*}$\frac{\kappa}{N}>1$}
\]  and
\[
 v=0  \quad \mbox{if\phantom{*}$\frac{\kappa}{N}<1$}
 . \]
\end{theo}
In particular, this implies that the positivity of the speed of a spider only depends on the
number of legs $N$ and not on the set $L$. Our technique is not fine enough to deal with the critical case $\kappa=N$. Nevertheless, we are inclined to believe that in this case, independently of the set $L$, the speed of the spider should be zero.

\section{Notations and auxiliary results}

We will denote by $K_1$, $K_2$, $\dots$ the ``important'' constants (those that can be used far away from the
place where they appear for the first time) and by $C_1$, $C_2$, $\dots$ the ``local'' ones (those that are used
only in a small neighbourhood of the place where they appear for the first time), restarting the numeration at
the beginning of each section in the latter case.

An important ingredient of our proofs is the analysis of the potential associated to the environment, which was
introduced by Sinai in \cite{Sinai}. The potential, denoted by $V=(V(x), x \in \Z)$ is a function of the
environment and is defined as follows:
\[
V(x) = \left\{
    \begin{array}{lll}
     \sum_{i=0}^{x-1}\ln\frac{\omega^-_i}{\omega^+_i}, & x>0,\\
     0, & x=0,\\
      \sum_{i=x+1}^{0}\ln\frac{\omega^+_i}{\omega^-_i}, & x<0.
    \end{array}
           \right.
\]

\subsection{Reversible measure of a spider}
\label{revmes} 
\begin{figure}
\begin{center}
\includegraphics[scale= 0.8]{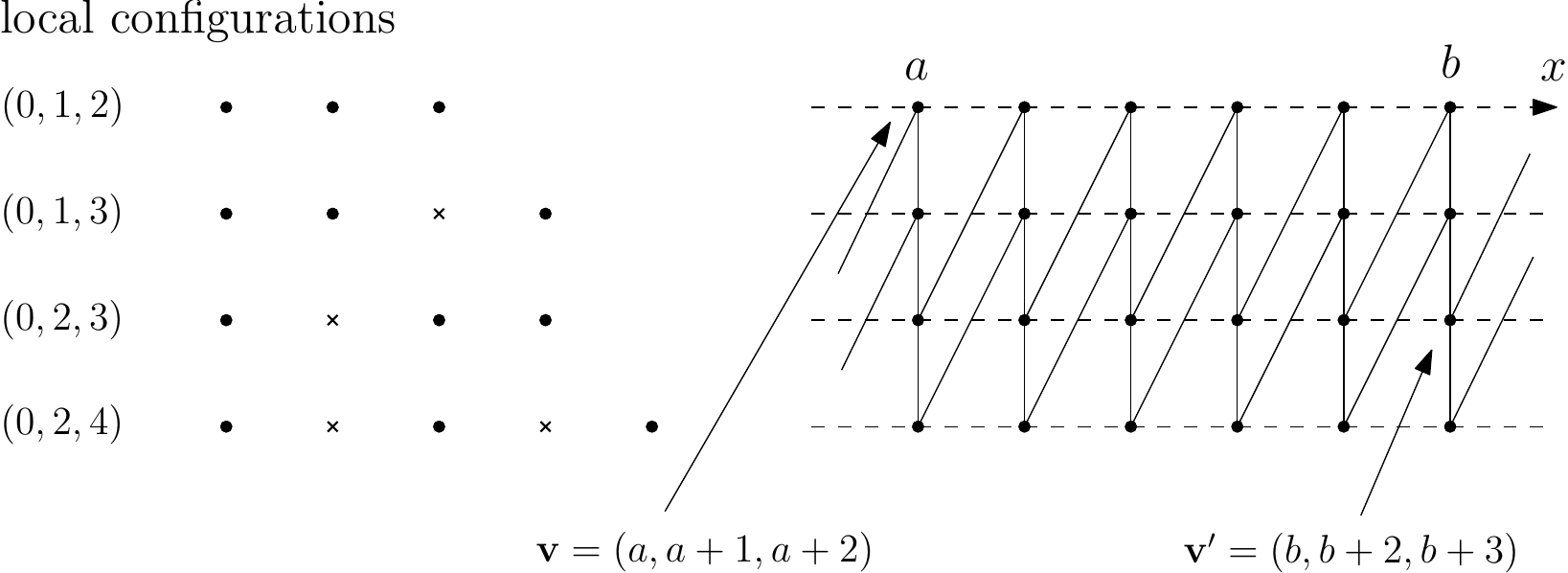}
\caption{Structure of the spider graph $\G$. The elements of $L$ are represented on the left.} \label{fig1}
\end{center}
\end{figure}

Let us first give an example to illustrate the construction of the spider graph $\G$. Consider a spider with 3 legs, that is $N=3$, and the following set $L$ of
restrictions:
 \[
 L=\{(0,1,2),(0,1,3),(0,2,3),(0,2,4)\}.
 \]
Figure \ref{fig1} shows the set of local configurations and a part of the spider graph $\G$. In the spider graph, the horizontal axis corresponds to the positions of the spider and the vertical axis to the local configurations. While the spider graph is deterministic the transition rates associated to each edge of
$\G$ depend on the realization $\omega$ of our random environment.

Now, given a couple $(N,L)$, consider  the continuous time Markov process $S=(S(t))_{t\geq 0}$ on the spider
graph $\G$. Observe that we use the same notation $S(t)_{t\geq 0}$ for two different processes: the spider on $\Z^N$ and the Markov process on $\G$. It should always be clear from the context to which of these we are
referring to.

Let   $\theta_x=e^{-V(x)}+e^{-V(x-1)}$. Note that $\theta_x$ is  the reversible measure at point $x$ for a single random walk on $\Z$. Then, the process $S$ is
$\P$-a.s.\ reversible with reversible measure
\begin{equation}
\label{revmesinf} \pi(\x)=\prod_{i=1}^{N} \theta_{x_i}
\end{equation}
for all $\x \in \V$. Using condition (iv) we obtain that
for all $x\in \Z$
\[
K_1e^{-V(x)}\leq \theta_{x} \leq K_2e^{-V(x)}
\]
for $K_1$ and $K_2$ two positive constants and
\[
 |V(x+1)-V(x)|\leq \ln \frac{1-\delta}{\delta}.
\]
Using these inequalities and the fact that $L$ is finite, we obtain that there exists two finite positive
constants $K_3$ and $K_4$ such that
\begin{equation}
\label{medrev} K_3e^{-NV(x_1)}\leq \pi(\x)\leq K_4e^{-NV(x_1)}
\end{equation}
for all $\x \in \V$. Now, let $I=[a,b]\cap \Z$ be a finite interval. Consider the graph $\G_I=(\V_I,
\Ed_I)\subset \G$ with
\begin{equation}
\label{graphvert}
\V_I=\bigcup_{x\in I}L_x
\end{equation}
\begin{equation}
\label{graphedge}
\Ed_I=\Big\{e=(\v,\w)\in \Ed\phantom{*} \mbox{such that}\phantom{*}\v,\w \in  \V_I \Big\}.
\end{equation}
Then, consider the process $\hat S$ which is the restriction of the process $S$ on the graph $\G_I$. As the graph
$\G_I$ is a subgraph of $\G$, the reversible measure (\ref{revmesinf}) is also reversible for the process $\hat
S$. Moreover as the graph $\G_I$ is finite we can normalize the reversible measure (\ref{revmesinf}) to obtain
the invariant probability measure $\hat \pi$ of $\hat S$,
\begin{equation}
\label{invmes} \hat{\pi}(\x)=\Big(\sum_{\x \in \V_I}\prod_{j=1}^{N} \theta_{x_j}\Big)^{-1}    \pi(\x)
\end{equation}
for all $\x \in \V_I$.

\subsection{Transience of the spider}
\label{transspid} Solomon showed in \cite{solom} that, under condition (iii), a single random walk is $\IP$-a.s.\
transient to $+\infty$. The following proposition shows that this is still the case for a spider.

\begin{figure}
\begin{center}
\includegraphics[scale= 0.9]{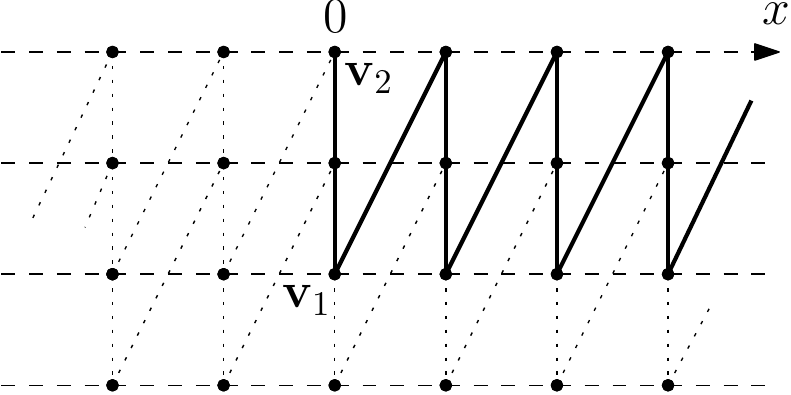}
\caption{Example of graph $\G'$ (in straight lines) for the spider of Figure~\ref{fig1}.} \label{fig4}
\end{center}
\end{figure}

\begin{prop}\label{prop:trans} Under the hypothesis (i)-(iv) a spider is always transient, that is,
\[
\lim_{t \to \infty} S_1(t)=\infty, \mbox{ $\IP$-a.s.}
\]
\end{prop}
{\it Proof.}\\
Consider the electrical network associated to $\G(\omega)$ by putting on each edge $e=(\x, \y)\in \Ed$ the
resistance $R_e=R_{\x,\y}=(q^S(\x,\y) \pi(\x))^{-1}$. By condition (ii), there exists in $L$ a vertex $\v_1$
which is linked to $L_1$ by some edge in $\Ed$. In the same way, there exists also a vertex $\v_2$ which is
linked to $L_{-1}$ by some edge in $\Ed$. By conditions (i) and (ii), we can choose a path
$\gamma_0$ from $\v_1$ to $\v_2$. As $\G$ is homogeneous, we can iterate this construction to all
the sets $L_x$, $x\geq 1$ and thus consider the linear sub-electrical network $\G'(\omega)=(\V',\Ed')$, see Figure~\ref{fig4}. If $I$ is an interval of $\N$, we define
\begin{equation*}
\V'(I)=\V_I\cap \V'
\end{equation*}
and
\begin{equation*}
\Ed'(I)=\Big\{e=(\v,\w)\in \Ed' \phantom{*} \mbox{such that} \phantom{*} \v,\w \in \V'(I) \Big \}.
\end{equation*}

Now, it is easy to compute the resistance $R_{\infty}=R_{\infty}(\omega)$ of $\G'(\omega)$ and to show that it
is finite. By definition,
\begin{equation*}
R_{\infty}:=\lim_{n\to \infty}R_n
\end{equation*}
where
\begin{equation*}
R_n=\sum_{e \in \Ed'([0,n])}R_{e}.
\end{equation*}
By condition (iv), we have
\begin{equation*}
R_n\leq \frac{1}{\delta}\sum_{\y \in \V'([0,n])} \pi^{-1}(\y).
\end{equation*}
Using inequality (\ref{medrev}), we obtain
\begin{equation}
\label{series1} R_n\leq \frac{K_3|L|}{\delta}\sum_{i=0}^n e^{NV(i)}.
\end{equation}
Let us first show that $\lim_{n\to \infty} R_n<\infty$ $\P$-a.s.
As $V(x)$, for $x\geq 0$, is a sum of bounded i.i.d. random variables, by the Strong Law of Large Numbers, we have
\begin{equation*}
\lim_{n\to \infty}\frac{V(n)}{n}=\E[\ln \rho_0]<0,\phantom{***}\mbox{$\P$-a.s.}
\end{equation*}
Now, take $\eps>0$ sufficiently small such that $(\E[\ln \rho_0]+\eps)<0$. 
This implies that $V(n)<(\E[\ln \rho_0]+\eps)n$ $\P$-a.s. Then, the  general term $e^{NV(n)}$ of (\ref{series1}) is dominated by $e^{N(\E[\ln \rho_0]+\eps)n}$ which is the general term of a convergent series. This shows that
\begin{equation*}
R_{\infty}=\lim_{n\to \infty} R_n<\infty,\phantom{**}\mbox{$\P$-a.s.}
\end{equation*}

As $\G'(\omega)$ is a sub-network of $\G(\omega)$ with $\P$-a.s.\ finite resistance, by the Rayleigh's
Monotonicity Law (see for example Doyle and Snell \cite{DS}) we deduce that the effective resistance of $\G(\omega)$ is
$\P$-a.s.\ finite, which implies that a spider on $\G(\omega)$ is transient for $\P$-almost all $\omega$.
\begin{flushright}
$\square$
\end{flushright}

\begin{rem}
In fact, Proposition~\ref{prop:trans} does hold in a more general context. Assume condition (iv) and let $(N,L)$ define a
spider. Then, one can show that the RWRE is recurrent iff the spider is recurrent. This follows from the fact that one can show
that the RWRE and the spider are roughly equivalent as electrical networks. 
We refer to~\cite{GMP} where these questions are discussed for a general spider.
\end{rem}
\subsection{Upper bound on the probability of confinement}
\label{upprobconfin} 
\begin{figure}
\begin{center}
\includegraphics[scale= 0.7]{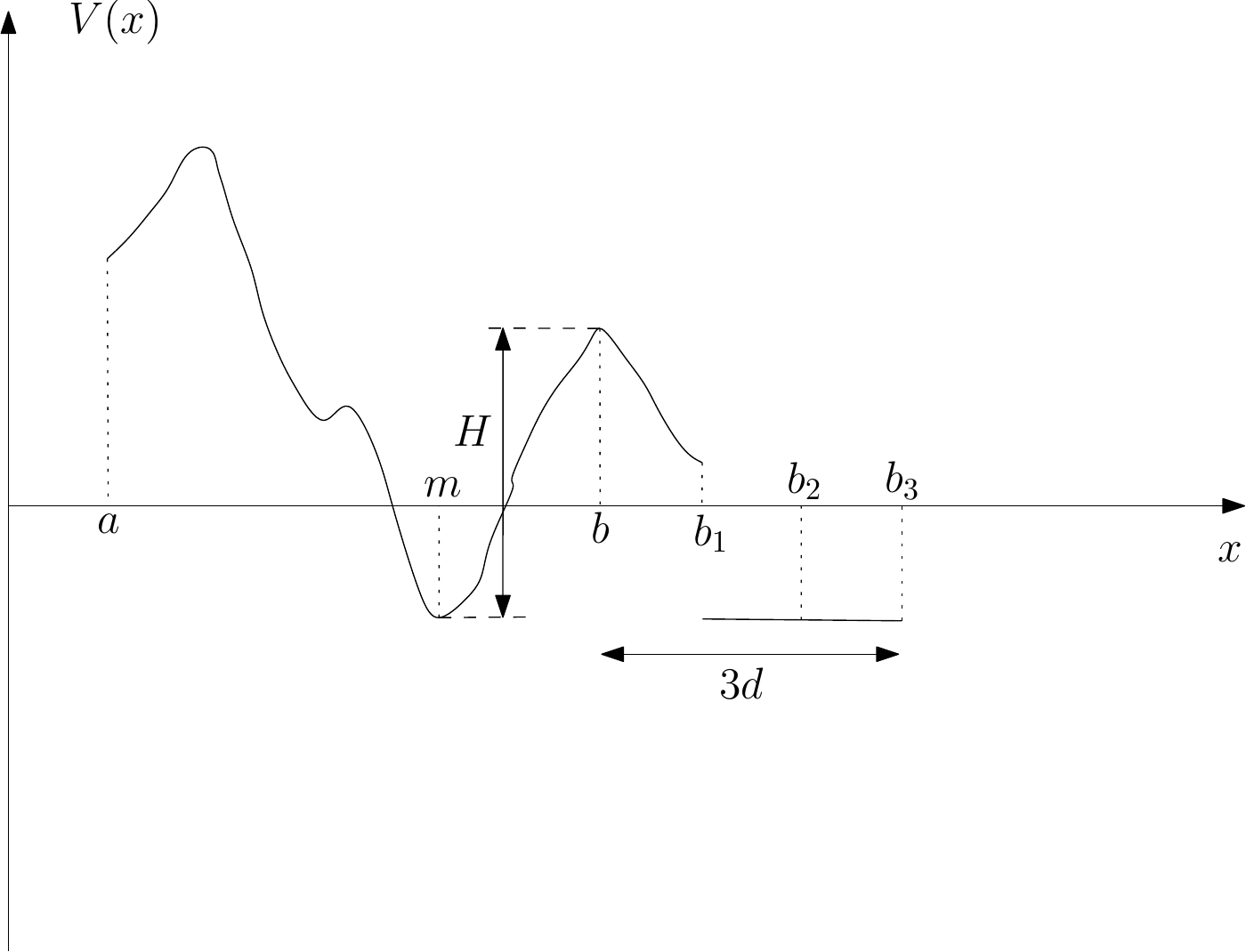}
\caption{Potential extention technique.} 
\label{fig3}
\end{center}
\end{figure}

In this section we want to deduce an upper estimate for the probability of confinement of a
spider on a finite interval. Fix a couple $(N,L)$ and let $I=[a,b]\cap \Z$, $a,b\in \Z$, be  a finite interval and
\[
\tau_{\{a,b\}}= \inf\{s>0: \mbox{$S_1(s)=a$ or $S_1(s)=b$}\}.
\]
We want to bound from above $\Po^{\x}[\tau_{\{a,b\}}>t]$ uniformly over all initial positions
$\x=(x_1,\dots,x_N)$ such that $a<x_1<b$. As $L$ is finite,
$d=\max_{\u,\v \in L}\|\u-\v \|_{\infty}$ is finite (where $\|\cdot\|_{\infty}$ is the usual $\infty$-norm in $\Z^N$). Let
$b_1=b+d$ and define $I_1=[a,b_1]\cap \Z$ and
\begin{equation}
\label{HACHE} H=\max_{x\in I_1}\Big(\max_{y\in[x,b_1]}V(y)-\min_{y\in[a,x)}V(y)\Big).
\end{equation}
Also, let
\begin{equation*}
m=\argmin_{x\in I_1}V(x).
\end{equation*}
We will show the following
\begin{prop}
\label{prop1} 
Let $[a,b]$ be a finite interval. We have
\begin{equation}
\label{upboundconf} \Po^{{\bf x}}[\tau_{\{a,b\}}>t]  \leq \exp \Big\{-\frac{t}{K_5(b-a)^5e^{NH}}\Big\}
\end{equation}
with $K_5$ a positive constant.
\end{prop}
{\it Proof.}\\
First we use the following trick: consider the interval $I_2:=[a,b_2]\supset I_1$, where $b_2=b+2d$ and an
interval $(b_2,b_3]$ such that $b_3-b_2=d$. On the subinterval $(b_1,b_3]$, we modify the environment such
that $V(x)=V(m)$ for every $x \in (b_1,b_3]$, see Figure \ref{fig3}.

Consider now the process $\hat{S}$ on the graph $\G_{I_2}$ (see (\ref{graphvert}) and (\ref{graphedge}) for the definition of $\G_{I_2}$) and define
\[\tau'=\inf\{s>0: \hat{S}(s)\in L_{b_2}\}.
\]
Since, $\Po^{\x}[\tau_{\{a,b\}}>t]\leq \Po^{\x}[\tau'>t]$, we focus from now on on finding
an upper bound for $\Po^{\x}[\tau'>t]$.

To this end, we construct a lower bound for the spectral gap $\lambda$ of the process $\hat{S}$ using Theorem 3.2.1 of Saloff-Coste~\cite{SC}. For all pairs of vertices $(\x, \y)$ of $\G_{I_2}$, choose a path
in $\Ed_{I_2}$ going from $\x$ to $\y$. We denote this path $\gamma(\x, \y)$ and let $\Gamma=\{\gamma(\x,
\y):(\x, \y)\in \V_{I_2} \times \V_{I_2}\}$. Then, the latter theorem states that $\lambda \geq 1/A$ where
\begin{equation}
\label{star}\phantom{**}A=\max_{e\in \Ed_{I_2}}\Big\{R_e\sum_{\x,\y\in \V_{I_2}: e\in \gamma(\x,\y)}|\gamma(\x, \y)|{\hat
\pi}(\x){\hat \pi}(\y)\Big\}
\end{equation}
and $R_e$ is the resistance of edge $e$ as defined in subsection \ref{transspid}. Now, let us define a set of paths
$\Gamma$ that will give a good lower bound for the spectral gap $\lambda$. We start by enumerating the elements
of the set $L$. The shift $\Theta$ induces the same enumeration on all the sets $L_x$ for $x\in \Z$. If
$\y \in L_x$ for some $x$, we will denote by $n(\y)$ the number associated to $\y$. Furthermore, let us fix two
local configurations $\r_1$ and $\r_2$ of $L$ such that the edge $e=(\Theta_{x}\r_1, \Theta_{x+1}\r_2)\in \Ed$
for all $x\in \Z$. Let $\x=(x_1,\dots,x_N)$ and $\y=(y_1,\dots,y_N)$ be two vertices of $\G_{I_2}$, we will now
choose a path $\gamma(\x,\y)$ as follows:
\begin{itemize}

\item if $\x$ and $\y$ are such that $x_1=y_1$ then consider the set of all paths that are contained in
$\Ed_{x_1}$ (see (\ref{graphedge}) for the definition of $\Ed_{x_1}$) which go from $\x$ to $\y$. Assume $n(\x)<n(\y)$. In this case, we choose the path $(\x, \x_1,\dots,
\x_N, \y)$ which minimizes the number $n(\x_1)\dots n(\x_N)$ in the following sense: $n(\x_1)\dots n(\x_{N_1})$ is smaller than $n(\x'_1)\dots n(\x'_{N_2})$ if $N_1<N_2$. If $N_1=N_2=N$, we use the lexicographical order to decide which is the smallest one, that is, $n(\x_1)\dots n(\x_{N})$ is smaller than $n(\x'_1)\dots n(\x'_{N})$ if there exists $k\leq N$ such that $n(\x_i)=n(\x'_i)$ for $i\leq k$ and $n(\x_k)<n(\x'_k)$. If $n(\x)>n(\y)$, define $\gamma(\x,\y)$ as the inverse path of $\gamma(\y,\x)$;

\item if $\x$ and $\y$ are such that $x_1\neq y_1$. Assume first that $x_1<y_1$. Observe that there exists an
element $\z$ of $\V_{x_1}$ such that $\z=\Theta_{x_1}\r_1$. Then, by the method above, we go from $\x$ to $\z$.
From $\z$, we go to $\x'\in \V_{x_1+1}$ such that $\x'=\Theta_{x_1+1}\r_2$. From now on, we iterate the
process to reach some $\z'$ such that $z'_1=y_1$. Finally, we again use
the method above to go from $\z'$ to $\y$. If $x_1>y_1$, define $\gamma(\x,\y)$ as the inverse path of $\gamma(\y,\x)$.
\end{itemize}
Thus, we have constructed the set $\Gamma$ we will use in the rest of this proof.

Now, let us find an upper bound of $A$ from~\eqref{star}. First let us define
\[
 A(e)=R_e\sum_{\x,\y \in \V_{I_2}: e\in \gamma(\x,\y)}|\gamma(\x,\y)|\hat{\pi}(\x)\hat{\pi}(\y)
\]
for all $e\in \Ed_{I_2}$. Let us find a uniform upper bound of $A(e)$ over all $e\in \Ed_{I_2}$. Let
$e=(\z,\w)$.  Using condition (iv), we obtain  that
\begin{equation*}
A(e)\leq \frac{1}{\delta \hat{\pi}(\z)}\sum_{\x,\y\in \V_{I_2}: e\in
\gamma(\x,\y)}|\gamma(\x,\y)|\hat{\pi}(\x)\hat{\pi}(\y).
\end{equation*}
Then, as $|\gamma(\x,\y)|$ is uniformly bounded by $|L| (b_2-a)$, using inequalities (\ref{medrev})
we obtain
\begin{equation*}
A(e)\leq C_3(b_2-a)D^{-1} e^{NV(z_1)}\sum_{\x,\y\in \V_{I_2}: e\in \gamma(\x,\y)}e^{-N(V(x_1)+V(y_1))},
\end{equation*}
where $D=\sum_{\x \in \V_{I_2}}\prod_{i=1}^{N} \theta_{x_i}$ and $C_3$ is a positive constant.

Now, using the rough upper bound 
\[
|\{\x,\y\in \V_{I_2}: e\in \gamma(\x,\y)\}|\leq (b_2-a+1)^2|L|^2,
\]
 by the
construction of $\Gamma$, we have
\begin{equation*}
A(e)\leq C_4(b_2-a)^3D^{-1}\exp\{NV(z_1)-N(\min_{x_1\leq z_1}V(x_1)+\min_{y_1\geq w_1}V(y_1))\}
\end{equation*}
with $C_4$ a positive constant. Now, observe that by (\ref{medrev})
\[D^{-1}\exp\{-N\min_{y_1\geq w_1}V(y_1)\}\leq \frac{1}{K_3}\]
and by definition of $H$ (see (\ref{HACHE})) 
\[\max_{z_1 \in I_2}[V(z_1)-\min_{x_1\leq z_1}V(x_1)]\leq H.\]
We obtain
\begin{equation*}
A(e)\leq C_5(b_2-a)^3e^{NH}.
\end{equation*}
By condition (iv) note that there exists a positive constant $C_6$ such that $(b_2-a)\leq C_6(b-a)$. Thus, we
obtain
\begin{equation*}
A=\max_{e\in \Ed_{I_2}}A(e)\leq C_7(b-a)^3e^{NH}.
\end{equation*}
and with Theorem 3.2.1 of~\cite{SC},
\begin{equation}
\label{spectgap} \lambda \geq \frac{1}{C_7(b-a)^3e^{NH}}.
\end{equation}
We are aiming now for a (uniform in $\x\in \V_{I_2}$) lower bound for
$\Po^{\x}[{\hat S}(s)\in L_{b_2}]$.
First, we recall the following fact: for $\x,\y\in
\V_{I_2}$ and $s>0$,
\begin{equation}
\label{conv-sp} \Big|\Po^{\x} [{\hat S}(s) = \y] - {\hat \pi}(\y)\Big| \leq \Big(\frac{{\hat \pi}(\y)}{{\hat
\pi}(\x)}\Big)^{1/2} \exp\{ -\lambda s \},
\end{equation}
see Corollary 2.1.5 in~\cite{SC}.  Furthermore, notice that
\begin{equation*}
\Po^{\x}[{\hat S}(s)\in L_{b_2}]\geq \Po^{\x}[{\hat S}(s)=\v]
\end{equation*}
for any $\v \in L_{b_2}$.\\
Then,  by inequalities (\ref{medrev}) and condition (iv)  for $\x$ such that $a<x_1<b$ and $\y$
such that $y_1=b_2$ we have
\[
\Big(\frac{{\hat \pi}(\y)}{{\hat \pi}(\x)}\Big)^{1/2}\leq
\Big(\frac{K_4}{K_3}\Big)^{\frac{1}{2}}e^{\frac{N}{2}(V(x_1)-V(b_2))}\leq e^{C_8(b-a)},
\]
where $C_8$ is a positive constant to be chosen later. Note that we can take $C_8$ arbitrary large.  Hence, using inequality (\ref{spectgap}) and  taking
\[
s= 2C_7C_8(b-a)^4e^{NH}
\]
we obtain
\begin{equation}
\Big(\frac{{\hat \pi}(\y)}{{\hat \pi}(\x)}\Big)^{1/2}\exp\{ -\lambda s \}\leq e^{-C_8(b-a)}.
\end{equation}
Since the potential is constant and equals to $V(m)$ on the interval $[b_1,b_3]$, we obtain
\begin{equation*}
{\hat \pi}(\v)\geq \frac{1}{2|L|(b_2-a)}\geq \frac{1}{2C_6|L|(b-a)}.
\end{equation*}
Suppose that $C_8$ is large enough so that
\[
 e^{-C_8(b-a)}\leq \frac{1}{4|L|C_6(b-a)}.
\]
Using (\ref{conv-sp}), we obtain
\begin{equation}
 \label{missing}
\Po^{\x}[{\hat S}(s)=\v]\geq  \frac{1}{4C_6|L|(b-a)}.
\end{equation}
Now, divide the time interval $[0,t]$ into $M:=\lfloor\frac{t}{s}\rfloor$ subintervals of length $s$. Using~\eqref{missing}   and
the Markov property we obtain
\begin{align*}
\Po^{\x}[\tau_{\{a,b\}}>t] &\leq  \Po^{\x}[\tau'>t] \nonumber\\
&\leq \Po^{\x}[{\hat S}(sj)\notin L_{b_2}, j=1,\dots,M] \nonumber\\
&\leq \Big(1-\frac{1}{4C_6|L|(b-a)}\Big)^M\nonumber\\
&\leq \exp \Big\{-\frac{M}{4C_6|L|(b-a)}\Big\} \nonumber\\
&\leq \exp \Big\{-\frac{t}{C_9(b-a)^5e^{NH}}\Big\}
\end{align*}
with $C_9$ a positive constant. \\
This concludes the proof of Proposition \ref{prop1}.
\begin{flushright}
$\square$
\end{flushright}

\subsection{Probability of escape in a given direction}

We also need the following result. For $y\in \Z$, let
\begin{equation}
\label{Pronf}
\tau_y=\inf \{s>0: S_1(s)=y\}.
\end{equation}
We can adapt Lemma 3.4 of  Comets and Popov~\cite{CP} in an elementary way to obtain the following upper bound for the probability of escape in a given direction.
\begin{prop}
\label{ESC} For some $K_6\in(0,\infty)$, we have for all $s>0$, $\x\in \V$, $y\in \Z$
\[\Po^{\x}[\tau_y<s]\leq K_6\int_0^{s+1}\Po^{\x}[S_1(u)=y]du.\]
\end{prop}

\section{Case $\kappa/N>1$}

This section is devoted to the proof of the positiveness of the speed of  a spider when $\kappa/N>1$.

\begin{figure}
\begin{center}
\includegraphics[scale= 0.6]{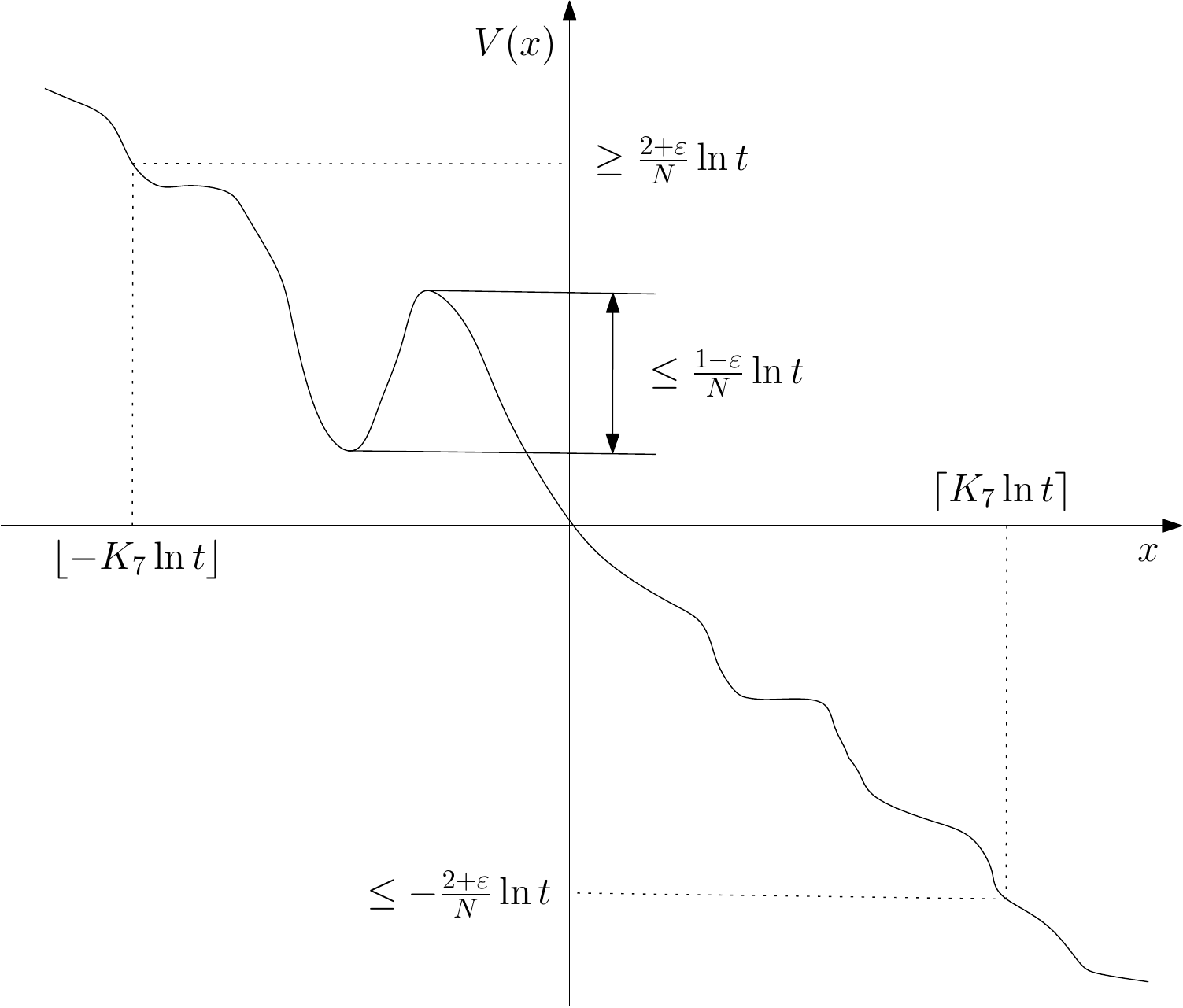}
\caption{On the definition of $\Lambda_t$.} 
\label{fig6}
\end{center}
\end{figure}

Fix a couple $(N,L)$ and let $\x_0\in L$ be an initial configuration of the spider. In order to simplify
notations, we will systematically omit the superscript $\x_0$ for the quenched and the annealed laws and
expectations. Remember that
\[
\T:= \inf \{s>0: S_1(s)>0 \phantom{*}\mbox{and}\phantom{*}S(s)=\Theta_{S_1(s)} \x_0\}.
\]
We will show that if $\frac{\kappa}{N}>1$ then $\IE[\T]<\infty$ which will imply by the Birkhoff's Ergodic
Theorem that $v>0$. First, for each $t>1$, we define the set of ``$t$-good'' environments.
\begin{df}
\label{def_goo_env} Fix $t>1$ and let $0<\eps<1$. Then, fix a finite absolute constant $K_7>0$ (i.e. $K_7$ does
not depend on $\omega$ and $t$). A realization of the potential $V$ is said to be $t$-good if we have
\begin{itemize}
\item $V(\lfloor -K_7\ln t \rfloor)\geq \frac{2+\eps}{N}\ln t$, \item $V(\lceil K_7\ln t \rceil)\leq
-\frac{2+\eps}{N}\ln t$, \item $\max_{i\in [\lfloor -K_7\ln t \rfloor, \lceil K_7\ln t \rceil]}\max_{j\geq
i}(V(j)-V(i))\leq \frac{1-\eps}{N}\ln t.$
\end{itemize}
We will call $\Lambda_t$ the set of $t$-good environments. See Figure \ref{fig6}.
\end{df}
The following decomposition is the key of our analysis.
\begin{align}
\label{decomp}
 \IP[\T>t]&=\int_{\Omega}\Po[\T>t]d\P(\omega)\nonumber\\
&\leq \sup_{\omega \in \Lambda_t}\Po[\T>t]+\P[\Lambda_t^c].
\end{align}
In the two following subsections we will show that both terms of the right-hand side of (\ref{decomp}) are
integrable in $t$ and thus $\IE[\T]<\infty$. We start with the term $\P[\Lambda_t^c]$.
\subsection{Upper bound on $\P[\Lambda_t^c]$}
By definition of $\Lambda_t$ we obtain that
\begin{align*}
\P[\Lambda_t^c] &\leq
\P\Big[V(\lfloor -K_7\ln t \rfloor)< \frac{2+\eps}{N}\ln t\Big] +\P\Big[V(\lceil K_7\ln t \rceil)> -\frac{2+\eps}{N}\ln t\Big]\nonumber\\
& \phantom{**}+\P\Big[\max_{i\in [\lfloor -K_7\ln t \rfloor,\lceil K_7\ln t \rceil]}\max_{j\geq i}(V(j)-V(i))>\frac{1-\eps}{N}\ln t\Big]\nonumber\\
&\leq 2\P\Big[V(\lceil K_7\ln t \rceil)> -\frac{2+\eps}{N}\ln t\Big]\nonumber\\
& \phantom{**}+\P\Big[\max_{i\in [0,2\lceil K_7\ln t \rceil]}\max_{j\geq i}(V(j)-V(i))>\frac{1-\eps}{N}\ln
t\Big].
\end{align*}
Let us define
\begin{equation*}
A_t=\Big\{V(\lceil K_7\ln t \rceil)> -\frac{2+\eps}{N}\ln t\Big\}
\end{equation*}
and
\begin{equation*}
B_t=\Big\{\max_{i\in [0,2\lceil K_7\ln t \rceil]}\max_{j\geq i}(V(j)-V(i))>\frac{1-\eps}{N}\ln t\Big\}.
\end{equation*}
Now, we will show that we can choose $K_7$ large enough such that  $\int_{0}^{\infty}\P[A_t] dt$ is finite.
Observe that as $\eps<1$ we have
\begin{align*}
\P[A_t]&= \P\Big[V(\lceil K_7\ln t \rceil)>-\frac{2+\eps}{N}\ln t\Big]\nonumber\\
&\leq \P\Big[V(\lceil K_7\ln t \rceil)>-\frac{3}{N}\ln t\Big]\nonumber\\
&\leq \P\Big[\frac{|V(\lceil K_7\ln t \rceil)-\E[V(1)]\lceil K_7\ln t \rceil|}{\lceil K_7\ln t \rceil}>a\Big]\nonumber\\
\end{align*}
for $a=-\frac{\E[V(1)]}{2}$ if $K_7>-\frac{6}{\E[V(1)]N}$. As $V(x)$, $x>0$, is a sum of bounded i.i.d.\
random variables, we can apply Cram\'er's Theorem to obtain that
\begin{equation*}
\P\Big[\frac{|V(\lceil K_7\ln t \rceil)-\E[V(1)]\lceil K_7\ln t \rceil|}{\lceil K_7\ln t \rceil}>a\Big]\leq
e^{-I(a)K_7\ln t}
\end{equation*}
with $I(\cdot)$ the large deviation function defined as
\[
I(x)=\sup_{l>0}[lx-\ln\E[lV(1)]].
\]
Taking $K_7>\frac{1}{I(a)}\vee-\frac{6}{\E[V(1)]N}$, we obtain that
\begin{equation}
\label{env1}
 \P[A_t]\leq e^{-C_1\ln t}= \frac{1}{t^{C_1}}
\end{equation}
with $C_1>1$. This shows that $\int_{0}^{\infty}\P[A_t] dt$ is finite.

Now, let us show that, with the choice of $K_7$ above, $\int_{0}^{\infty} \P[B_t] dt$ is finite too. We have
\begin{align*}
 \P[B_t]&=\P \Big[\max_{i\in [0,2\lceil K_7\ln t \rceil]}\max_{j\geq i}(V(j)-V(i))>\frac{1-\eps}{N}\ln t\Big] \nonumber\\
&\leq \sum_{i=0}^{2\lceil K_7\ln t \rceil}\P \Big[\max_{j\geq i}(V(j)-V(i))>\frac{1-\eps}{N}\ln t\Big].
\end{align*}
The estimate (2.7) in \cite{PGF} yields
\begin{equation}
\label{env2}
 \P[B_t]\leq C_2\frac{\lceil K_7\ln t \rceil}{t^{\frac{\kappa}{N}(1-\eps)}}
\end{equation}
with $C_2$ a positive finite constant. As $\frac{\kappa}{N}>1$, we can choose $\eps$ sufficiently small such
that $\frac{\kappa}{N}(1-\eps)>1$, this shows the integrability of $\P[B_t]$.
\subsection{Upper bound on $\sup_{\omega \in \Lambda_t}\Po[\T>t]$}
Let us denote $x_l=\lfloor -K_7\ln t \rfloor$ and $x_r=\lceil K_7\ln t \rceil$. Recall  that the initial configuration of the spider is $\x_0\in L$. We use the following
decomposition
\begin{align}
\label{decomp2}
 \Po[\T>t]
 &= \Po\Big[\T>t, \tau_{x_r}>\frac{t}{2}\Big]   +    \Po\Big[\T>t, \tau_{x_r}\leq \frac{t}{2}\Big] \nonumber\\
 &\leq \Po\Big[\tau_{x_r}>\frac{t}{2}\Big]   +    \Po\Big[\T>t, \tau_{x_r}\leq \frac{t}{2}\Big].
\end{align}

\subsubsection{Upper bound on $\Po\Big[\tau_{x_r}>\frac{t}{2}\Big] $}
 We write
\begin{align}
\label{decomp3}
 \Po\Big[\tau_{x_r}>\frac{t}{2}\Big]
 &=  \Po\Big[\tau_{x_r}>\frac{t}{2},\tau_{x_r}>\tau_{x_l}\Big]   +    \Po\Big[\tau_{x_r}>\frac{t}{2},\tau_{x_r}< \tau_{x_l}\Big] \nonumber\\
 &\leq   \Po[\tau_{x_r}>\tau_{x_l} ] +    \Po\Big[ \tau_{\{x_l,x_r\}}>\frac{t}{2}\Big].
 \end{align}
 Let us first treat the second term of the right-hand side of (\ref{decomp3}). Observe that on the interval $[x_l,x_r+d]$ (where $d$ is from subsection \ref{upprobconfin}), we have for $\omega \in \Lambda_t$, $e^{NH}\leq C_1t^{1-\eps}$ with $C_1$ a positive constant.
Using Proposition \ref{prop1} we obtain immediately that
\begin{equation}
\label{UP1} \Po\Big[\tau_{\{x_l,x_r\}}>\frac{t}{2}\Big]\leq \exp \Big\{-\frac{t^{\eps}}{C_2(\ln t)^5}\Big\}
\end{equation}
with $C_2$ a positive constant.

For the first term of the right-hand side of (\ref{decomp3}) let us write
\begin{align}
\label{decomp4}
\Po[\tau_{x_r}>\tau_{x_l}]&=\Po\Big[\tau_{x_r}>\tau_{x_l},\tau_{\{x_l,x_r\}}>\frac{t}{2}\Big]+\Po\Big[\tau_{x_r}>\tau_{x_l},\tau_{\{x_l,x_r\}}\leq \frac{t}{2}\Big] \nonumber\\
&\leq \Po\Big[\tau_{\{x_l,x_r\}}>\frac{t}{2}\Big]+\Po[\tau_{x_l} < t].
\end{align}
We can bound from above the first term of the right-hand side of (\ref{decomp4}) using (\ref{UP1}).

In order to bound from above the second term (\ref{decomp4}), we use Proposition \ref{ESC} and (\ref{medrev}) to obtain
\begin{align*}
\Po[\tau_{x_l} < t]
&\leq K_6 \int_0^{t+1} \Po^{\x_0}[S_1(u)=x_l]du \nonumber\\
& = K_6 \int_0^{t+1} \sum_{\y \in L_{x_l}} \Po^{\x_0}[S(u)=\y]du \nonumber\\
& = K_6 \int_0^{t+1} \sum_{\y \in L_{x_l}} \frac{\pi(\y)}{\pi(\x_0)}\Po^{\y}[S(u)=\x_0]du \nonumber\\
& \leq K_6|L|(t+1)\frac{\pi(\y)}{\pi(\x_0)}\nonumber\\
& \leq C_3 t e^{-NV(x_l)}
\end{align*}
with $C_3$ a positive constant.\\
For $\omega \in \Lambda_t$, we have that $V(x_l)> \frac{2+\eps}{N}\ln t$. Hence,
\begin{equation}
\label{UP2} \Po[\tau_{x_l} < t]\leq \frac{C_4}{t^{1+\eps}}
\end{equation}
with $C_4$ a positive constant. Eventually, by (\ref{decomp3}), (\ref{UP1}), (\ref{decomp4}) and (\ref{UP2}), we
obtain
\begin{equation}
\label{TY1}
 \Po\Big[\tau_{x_r}>\frac{t}{2}\Big] \leq 2\exp \Big\{-\frac{t^{\eps}}{C_2(\ln t)^5}\Big\}+\frac{C_4}{t^{1+\eps}}
\end{equation}
for $\omega \in \Lambda_t$.

\subsubsection{Upper bound on $\Po\Big[\T>t, \tau_{x_r}\leq \frac{t}{2}\Big]$}

Let $\mathfrak{F}_{x_r}$ be the $\sigma$-field generated by the process $S$ up to the stopping time
$\tau_{x_r}$. Using the Markov property, we obtain
\begin{align*}
 \Po\Big[\T>t, \tau_{x_r}\leq \frac{t}{2}\Big]
 & =\Eo\Big[1_{\{\tau_{x_r}\leq t/2\}}\Po[\T>t\mid  \mathfrak{F}_{x_r}] \Big]\nonumber\\
 &\leq \Eo\Big[1_{\{\tau_{x_r}\leq t/2\}}\Po^{S(\tau_{x_r})}\Big[\T>\frac{t}{2}\Big] \Big]\nonumber\\
 &\leq \Po\Big[\tau_{x_r}\leq \frac{t}{2}\Big] \times \max_{\y \in L_{x_r}}\Po^{\y}\Big[\T>\frac{t}{2}\Big].
\nonumber\\
 &\leq \max_{\y \in L_{x_r}}\Po^{\y}\Big[\T>\frac{t}{2}\Big].
\end{align*}
The next step is to bound uniformly in $\y$ the quantity $\Po^{\y}\Big[\T>\frac{t}{2}\Big]$ for $\y \in
L_{x_r}$. We use the following decomposition
\begin{align}
\label{DURT} \Po^{\y}\Big[\T>\frac{t}{2}\Big]
&= \Po^{\y}\Big[\T>\frac{t}{2},\tau_0< t\Big]+\Po^{\y}\Big[\T>\frac{t}{2},\tau_0\geq t\Big]\nonumber\\
&\leq \Po^{\y}[\tau_0< t]+\Po^{\y}\Big[\T>\frac{t}{2},\tau_0\geq t\Big].
\end{align}

To bound from above the first term of the right-hand side of (\ref{DURT}), we use Proposition \ref{ESC} to
obtain
\begin{equation}
\label{Tyff}
\Po^{\y}[\tau_0< t]\leq \frac{C_4}{t^{1+\eps}}
\end{equation}
for $\omega \in \Lambda_t$.

For the second term of the right-hand side of (\ref{DURT}) we start by defining
\[\T'=\inf\{s>0: S(s)= \Theta_{S_1(s)}\x_0\}.
\]
Then, let $\Upsilon$ be the number of movements it takes the spider to be in local configuration $\x_0$ for the first time.\\ Formally, if $\Xi=(\Xi(n))_{n\geq 0}=(\Xi_1(n),\dots,\Xi_N(n))_{n\geq 0}$ is the
jump chain (recall that the jump chain of a jump Markov process is the sequence of states
visited by the Markov process) associated to the jump Markov process $S$ we have
\[ \Upsilon=\min\{n\geq1: \Xi(n)=\Theta_{\Xi_1(n)}\x_0\}.\]
 Observe that by condition (iv) and the facts that the process is irreducible and $L$ is finite, $\Upsilon<\infty$ $\IP$-a.s.
 Now, let $(T_i)_{i\geq 1}$ be the sequence of jump times of the process $S$. Observe that
 \begin{align}
 \label{URT0}
\lefteqn{\Po^{\y}[\exists s\in [0,2|L|]: S(s)= \Theta_{S_1(s)}\x_0]}\phantom{***}\nonumber\\
&= \Po^{\y}[T_1+T_2+\dots+T_{\Upsilon}\leq 2|L|]\nonumber\\
&\geq \Po^{\y}[T_1+T_2+\dots+T_{\Upsilon}\leq 2|L|, \Upsilon \leq |L|]\nonumber\\
&\geq  \Po^{\y}[T_1+T_2+\dots+T_{|L|}\leq 2|L|, \Upsilon\leq |L|]\nonumber\\
&=  \Po^{\y}[T_1+T_2+\dots+T_{|L|}\leq 2|L|]\Po^{\y}[ \Upsilon\leq |L|].
 \end{align}
By condition (iv), there exists $\eta>0$ such that
\begin{equation}
\label{URT1} \Po^{\y}[\Upsilon\leq |L|]\ge \eta
\end{equation}
for all $\y$.
Using Markov's inequality and condition (iv) we have
\begin{align}
\label{URT2}
 \Po^{\y}[T_1+T_2+\dots+T_{|L|}\leq 2|L|]
 &\geq 1-\frac{|L|}{2|L|}\nonumber\\
 &\geq 1-\frac{1}{2}\nonumber\\
 &\geq \frac{1}{2}.
 \end{align}
Therefore, using (\ref{URT1}), (\ref{URT2}) and  (\ref{URT0}) we obtain
\begin{equation*}
\max_{\y\in L_{x_r}}\Po^{\y}[\exists s\in [0,2|L|]: S(s)= \Theta_{S_1(s)}\x_0]\geq \frac{\eta}{2}.
\end{equation*}
The next step is to divide the interval $[0,\frac{t}{2}]$ into $\lfloor \frac{t}{4|L|} \rfloor$ intervals of size $2|L|$ and
observe that by the Markov property,
\begin{equation*}
 \Po^{\y}\Big[\T'>\frac{t}{2}\Big]  \leq \Big(1-\frac{\eta}{2}\Big)^{\lfloor \frac{t}{4|L|} \rfloor}
\end{equation*}
for all $\y \in L_{x_r}$. As $\{S(0)=\y,\T>t/2,\tau_0> t \} \subset \{S(0)=\y,\T'>\frac{t}{2}\}$ for all $\y \in L_{x_r}$ and $t$ sufficiently large, we obtain
\begin{equation}
\label{TY2}
 \Po^{\y}\Big[\T>\frac{t}{2},\tau_0> t\Big]  \leq \Big(1-\frac{\eta}{2}\Big)^{\lfloor \frac{t}{4|L|} \rfloor}
\end{equation}
for all $\y \in L_{x_r}$.

To sum up, by   (\ref{TY1}), (\ref{Tyff}) and  (\ref{TY2}) we obtain
\begin{equation}
\label{supom} \sup_{\omega \in \Lambda_t}\Po[\T>t]\leq 2\exp \Big\{-\frac{t^{\eps}}{C_2(\ln
t)^5}\Big\}+2\frac{C_4}{t^{1+\eps}}+\Big(1-\frac{\eta}{2}\Big)^{\lfloor \frac{t}{4|L|} \rfloor}.
\end{equation}
This shows that $\int_0^{\infty}\sup_{\omega \in \Lambda_t}\Po[\T>t]dt$ is finite.

\subsection{Positiveness of the speed}
In this subsection we show that the speed of the spider is positive $\IP$-a.s.\ if $\kappa/N>1$. Let $\zeta_0=0$
and
\[\zeta_n=\inf\{j>\zeta_{n-1}, \Xi_1(j)> \Xi_1(\zeta_{n-1}) \phantom{*}\mbox{and}\phantom{*} \Xi(j)=\Theta_{\Xi_1(j)}\x_0\}\]
for $n \geq 1$.\\
Since the sequence $(\zeta_{n+1}-\zeta_{n})_{n\geq 0}$ is ergodic under the annealed measure $\IP$, we can apply
the Birkhoff's Ergodic Theorem  to obtain that
\begin{equation}
\label{BET}
 \lim_{n \rightarrow \infty} \frac{\zeta_n}{n}=\IE[\zeta_1]=\IE[\T]<\infty,\phantom{*}\mbox{$\IP$-a.s.}
\end{equation}
where the last equality will be shown below.

Now, take $\zeta_n\leq m< \zeta_{n+1}$, we obtain
\begin{equation*}
\Xi_1(\zeta_n)-(\zeta_{n+1}-\zeta_n)\leq \Xi_1(m)<\Xi_1(\zeta_n)+(\zeta_{n+1}-\zeta_n)
\end{equation*}
 which implies
\begin{equation}
\label{RTY} \frac{\Xi_1(\zeta_n)-(\zeta_{n+1}-\zeta_n)}{n}\frac{n}{\zeta_{n+1}}\leq
\frac{\Xi_1(m)}{m}<\frac{\Xi_1(\zeta_n)+(\zeta_{n+1}-\zeta_n)}{n}\frac{n}{\zeta_n}.
\end{equation}
Observe that the sequence $(\Xi_1(\zeta_{n+1})-\Xi_1(\zeta_n))_{n\geq 0}$ is also ergodic under $\IP$. Therefore
we can apply the Birkhoff's Ergodic Theorem to obtain
\begin{equation*}
 \lim_{m \to \infty} \frac{\Xi_1(\zeta_m)}{m}=\IE[\Xi_1(\zeta_1)]=\IE[S_1(\T)]>0,\phantom{*}\mbox{$\IP$-a.s.}
\end{equation*}
where the last equality follows from the fact that $\Xi_1(\zeta_1)=S_1(\T)$.
 Now, let $m\to \infty$ in (\ref{RTY}). As the spider is transient to the right, we have also $n\to \infty$ $\IP$-a.s. Thus, we can deduce that
\begin{equation}
\label{EMBspeed} \lim_{m \to \infty}
\frac{\Xi_1(m)}{m}=\frac{\IE[S_1(\T)]}{\IE[\T]}>0,\phantom{**}\mbox{$\IP$-a.s.}
\end{equation}
\medskip

The result (\ref{EMBspeed}), obtained for the embedded Markov chain, transfers to continuous time. Indeed, there
exists a family $(e_i)_{i\geq1}$ of exponential random variables of parameter 1, such that the $n$th jump of the
continuous  time random process $S$ occurs at time $\sum_{i=1}^{n}e_i$. These random variables are independent
of the environment and the discrete-time random walk. It follows that we can write $\T=\sum_{i=1}^{\zeta_1}e_i$.
Hence, $\IE[\T]=\IE[\zeta_1]E[e_1]=\IE[\zeta_1]$.

Let us denote by $R_n$ the time of the $n$th jump of $S$. Then, take $R_n\leq t <R_{n+1}$, we obtain
\begin{equation*}
\frac{\Xi_1(n)}{R_{n+1}}\leq \frac{S_1(t)}{t}< \frac{\Xi_1(n)}{R_{n}}
\end{equation*}
and consequently
\begin{equation}
\label{GEN} \frac{\Xi_1(n)}{n}\frac{n}{R_{n+1}}\leq \frac{S_1(t)}{t}< \frac{\Xi_1(n)}{n}\frac{n}{R_{n}}.
\end{equation}
Eventually, taking the limit $n\to \infty$ in inequality (\ref{GEN}), using (\ref{EMBspeed}) and the law of large
numbers for the sequence $R_1$, $R_2-R_1$, $R_3-R_2$, $\dots$, we obtain
\begin{equation*}
v=\lim_{t \to \infty}\frac{S_1(t)}{t}= \frac{\IE[S_1(\T)]}{\IE[\T]}>0 \phantom{*}\mbox{$\IP$-a.s.}
\end{equation*}

\section{Case $\kappa/N<1$}
In this last section, we show that if $\kappa/N<1$, the speed of a spider is null. 
First, we need to introduce
some notations. Following Fribergh, Gantert and Popov~\cite{PGF} we will define the valleys of the potential $V$
as follows. For $t>1$, we define by induction the environment dependent sequence $(J_i(t))_{i\geq 1}$ by
\begin{equation*}
J_0(t)=0,
\end{equation*}
\begin{align*}
J_{i+1}(t) =\min\{j: \; j\geq J_i(t), & V(J_i(t))-\min_{l\in[J_i(t),j]}V(l)\geq \frac{3}{1\wedge \kappa}\ln t,\\
& V(j)=\max_{l\geq j} V(l)\}.
\end{align*}
In the following the dependence on $t$ will be frequently omitted to ease the notations. The portion of the
environment $[J_i,J_{i+1})$ is called the $i$th valley. In \cite{PGF}, it is shown that for $t$ large enough the
valleys are descending in the sense that $V(J_{i+1})<V(J_i)$ for all $i\geq 0$. Then, we define the depth of the
$i$th valley as follows
\begin{equation*}
H_i=\max_{J_i(t)\leq j<l<J_{i+1}(t))} (V(l)-V(j)).
\end{equation*}
\medskip

\noindent
\begin{figure}
\begin{center}
\includegraphics[scale= 0.6]{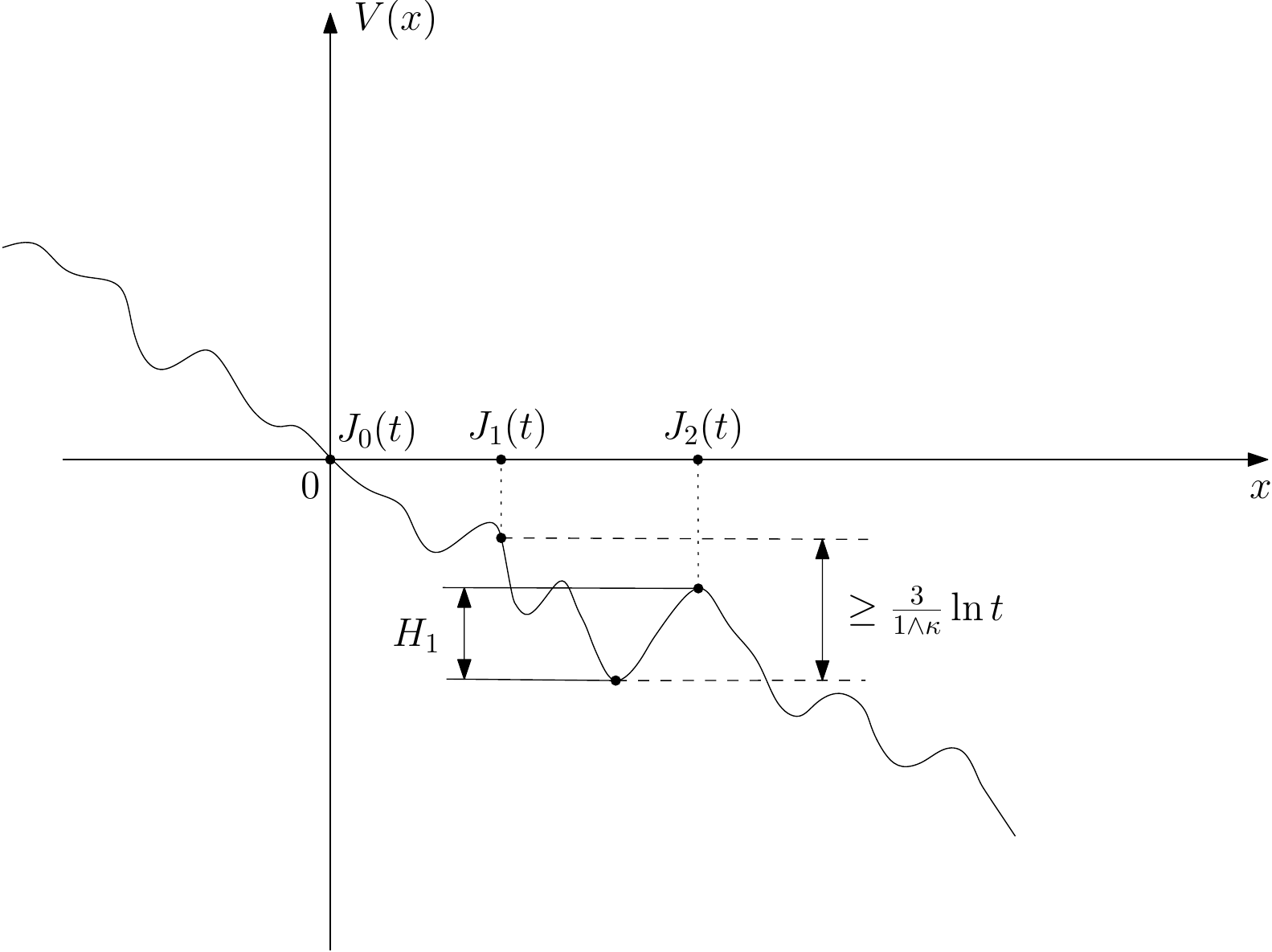}
\caption{On the definition of the valleys.}
 \label{fig2}
\end{center}
\end{figure}
Let us denote
\begin{equation*}
\L_t(m,m')=\{i\geq 1: [J_i,J_{i+1})\cap[\lfloor m \rfloor, \lfloor m' \rfloor) \neq \emptyset \}.
\end{equation*}

We define $\nu_0:=\frac{\kappa}{N}$ and consider $\nu$ such that $\nu_0<\nu<1$. Then,  take
$\eps=\eps(t)=\frac{4\ln\ln t}{\ln t}$ and define
\[
\mathcal{M}=\Big\{i\in \L_t(0,t^{\nu}): H_i\geq \frac{1-\eps}{\kappa}\ln t \Big\},
\]
\[ \Psi_t=\Big\{\omega:  |\mathcal{M}|\geq \frac{1}{3}t^{\nu-\nu_0}\Big\}.
\]
By Lemma 3.5 of \cite{PGF}, on each subinterval of length $t^{\nu_0}$, we find a valley of depth at least
$\frac{1-\eps}{\kappa}\ln t $ with probability at least $1/2$ for sufficiently large $t$. Since the interval
$[0,t^{\nu}]$ contains $t^{\nu-\nu_0}$ such intervals, we have
\begin{equation*}
\P[ \Psi_t]\geq 1-\exp(-C_1t^{\nu-\nu_0}).
\end{equation*}
For $i\in \mathcal{M}$, using the notation defined in (\ref{Pronf}), define $\sigma_i=\tau_{J_{i+1}+1}-\tau_{J_i+1}$ and let
\[s_0=\frac{t}{4\gamma_2(\ln t)^4}.
\]
Then, by Lemma 3.4 of \cite{CP} and the fact that $\kappa/N<1$, for any $i\in \mathcal{M}$,
\begin{align}
\label{FLUU} \Po^{\x}[\sigma_i<s_0]
&\leq 2\gamma_2s_0\exp \Big(- \frac{N}{\kappa}(\ln t -4\ln\ln t)  \Big) \nonumber\\
&\leq 2\gamma_2s_0\exp \Big(-\ln t +4\ln\ln t \Big) \nonumber\\
&= 2\gamma_2 s_0 t^{-1}(\ln t)^4\nonumber\\
&=\frac{1}{2}.
\end{align}
uniformly in $\x$, for sufficiently large $t$.

Define the family of random variables $\zeta_i= {\bf 1}\{\sigma_i < s_0\}$, $i\in \mathcal{M}$. Observe that by
the Markov property and (\ref{FLUU}), the sequence $(\zeta_i)_i$ is stochastically dominated by a sequence of
independent Bernoulli$\{0,1\}$ random variables $(\eta_i)_i$ of parameter $1/2$. Moreover, for $t$ large enough,
observe that we have
\[ \frac{1}{3}s_0 \frac{1}{3}t^{\nu-\nu_0}=\frac{1}{36\gamma_2(\ln t)^4}t^{1+\nu-\nu_0}>t.
\]
With the same notation as in (\ref{Pronf}), since $|\mathcal{M}| \geq \frac{1}{3}t^{\nu-\nu_0}$ for $\omega \in \Psi_t$ and since the sequence
$(\eta_i)_{i\geq 1}$ is i.i.d., we see using Cram\'er's theorem that for $t$ large enough
\begin{align*}
\Po[\tau_{\lfloor t^{\nu} \rfloor}<t]
&\leq \Po\Big[ \sum_{i\in \mathcal{M}}\zeta_i>\frac{2}{3}|\mathcal{M}|\Big]\nonumber\\
&\leq P\Big[ \sum_{i\in\mathcal{M}}\eta_i>\frac{2}{3} |\mathcal{M}|\Big]\nonumber\\
&\leq \exp \Big(-C_2t^{\nu-\nu_0}\Big).
\end{align*}
From this last inequality, we immediately conclude that the speed of the spider is null $\IP$-a.s.

\section*{Acknowledgements}
C.G. is grateful to Fapesp (grant  2009/51139-3) for financial support. S.M. thanks DFG (project MU 2868/1--1)
and Fapesp (grant 2009/08665-6) for financial support. S.P. and M.V. are grateful 
CNPq (grants 300886/2008--0, 472431/2009--9, 304561/2006--1) for financial support.  S.M., S.P., and M.V. thank
CAPES/DAAD (Probral) for support.

\end{document}